\documentclass[12pt]{amsart}
\usepackage[dvips]{graphicx}
\usepackage{amssymb}
\usepackage{amsmath}
\usepackage{latexsym}

\theoremstyle{plain}
\newtheorem{thm}{Theorem}[subsection]

\newtheorem{lemma}[thm]{Lemma}

\theoremstyle{definition}
\newtheorem{definition}[thm]{Definition}

\theoremstyle{remark}
\newtheorem{example}[thm]{Example}
\newtheorem{rem}[thm]{Remark}

\begin{document}

\title [statistically- and 2n-convergent series]
{Riemann Rearrangement Theorem for some types of convergence}

\author[Y.~Dybskiy]{Yuriy Dybskiy}
\address[Y.~Dybskiy]{Department of Mathematics \newline
Kharkov National University \newline
4, Svobody sq. Kharkov, UKRAINE}
\email{dybskiy@gmail.com}

\author[K.~Slutsky]{Konstantin Slutsky}
\address[K.~Slutsky]{Department of Mathematics \newline
Kharkov National University \newline
4, Svobody sq. Kharkov, UKRAINE}
\email{kslutsky@gmail.com}

\date{December 22, 2006}

\renewcommand{\theequation}{\thesubsection.\arabic{equation}}

\newcommand{\eqdef}{\stackrel{\mathrm{def}}{=}}
\newcommand{\eqst}{\stackrel{\mathrm{st.}}{=}}

\newcommand{\fseries}{\sum_\mathfrak{F}x_k}
\newcommand{\fsr}{SR_{\,2}(\sum x_k)}
\newcommand{\fsrs}{SR_{\,2}}
\newcommand{\R}{\mathbb{R}}
\newcommand{\N}{\mathbb{N}}
\newcommand{\Z}{\mathbb{Z}}
\newcommand{\Sg}{\mathfrak{S}}
\newcommand{\e}{\varepsilon}
\newcommand{\D}{\Delta}
\newcommand{\de}{\delta}

\newcommand{\stsr}{SR_{\,st.}(\sum x_k)}
\newcommand{\stsrs}{SR_{\,st.}}
\newcommand{\stat}{\stackrel{\mathrm{stat.}}{\longrightarrow}}

\begin{abstract}
We reexamine the Riemann Rearrangement Theorem for different types
of convergence. We consider series convergence with respect to a
filter. We describe the Sum Range (SR) of a series along the
2n-filter and for statistically convergent series.
\end{abstract}
\maketitle

\section{Introduction}

The classical Riemann Rearrangement Theorem says that the commutative law
is no longer true for infinite sums. To be more precise it says the
following:
\begin{thm}[Riemann Rearrangement Theorem]
Let $\sum_{k=1}^\infty x_k$ be a conditionally convergent series of
real numbers. Then:
\begin{enumerate}
  \item for any $s \in \R$ one can find a permutation $\pi$ such
  that \\ $\sum_{k=1}^\infty x_{\pi(k)} = s$;
  \item one can find permutation $\sigma$ such that $\sum_{k=1}^\infty x_{\sigma(k)} = \infty$;
  \item one can find permutation $\sigma$ such that $\sum_{k=1}^\infty x_{\sigma(k)} = -\infty$.
\end{enumerate}
\end{thm}

In the Riemann Rearrangement Theorem one considers the ordinary
convergence of series. It looks natural to consider in this setting
some weaker types of convergence. Interesting results in this
direction are proved in \cite{ces} and \cite{mat} where
generalizations of the Riemann theorem for Cesaro summation and
other matrix summation methods were obtained. These generalizations
are much more complicated than the original Riemann theorem, and
even the statements strongly differ from the classical one: say for
Cesaro summation it is possible that the set of sums under all
permutations of summands forms an arithmetic sequence. V.~Kadets
posed to us the problem what effects appear if the ordinary
convergence in the Riemann theorem statement is substituted by
convergence with respect to a filter. In this paper we are doing the
first two steps in this direction, considering statistical
convergence and convergence of subsequence $\sum_{k=1}^{2n} x_k$ of
partial sums.

\section{Statistical convergence}
\subsection{Introduction} Statistical convergence is a
generalization of the usual notion of convergence that parallels
the usual theory of convergence. While statistical convergence has
become an active area of research under the name of statistical
convergence only recently, it has appeared in the literature in a
variety of guises since the beginning of twentieth  century.
Statistical convergence has been discussed in number theory,
trigonometric series and summability theory. A relation between
statistical convergence and Banach space theory as well as a list
of references one can find in \cite{cgk2000}. The aim of this
chapter is to generalize the Riemann's Rearrangement Theorem to
the case of the statistical convergence.

The object that is going to be investigated is $\stsr$ and the
sequence of definitions below leads to it.

\begin{definition}
$A \subset \N$ is said to be {\it negligible} if
$$\lim_{n\rightarrow \infty} \frac{|A\cap\{1, \ldots , n\}|}{n} = 0$$
\end{definition}

\begin{definition}
The sequence $\{s_n\}_{n=1}^\infty$ {\it statistically converges} to
$s$ $(s_n \stat s)$ if for every $\e> 0$ the set $\{n:|s_n - s|> \e\}$ is
negligible.
\end{definition}

\begin{definition}
Series $\sum x_k$ is said to be {\it convergent statistically} to~$s$
if the sequence $s_n = \sum_1^n x_k$ of partial sums converges
statistically to $s$ (short notation is $\sum x_k\eqst s$).
\end{definition}

\begin{definition}
Point $s$ belongs to the {\it statistical Sum Range} of the series
$\sum x_k$ if there exists a permutation $\pi$ such that
$\sum_{k=1}^n x_{\pi(k)} \stat s$. The set of all such points is
called the {\it statistical Sum Range} of the series and is
denoted by $\stsr$.
\end{definition}

We will use also the following definition from \cite{sib}.
\begin{definition}
A point is said to be a {\it limit point} for the series $\sum x_k$
if it is the limit point of some subsequence of the sequence of
partial sums of some rearrangement of the series. The set of all
such points, called the {\it limit-point range} of the series, will
be denoted by LPR$(\sum x_k)$
\end{definition}

It is easy to see that LPR$(\sum x_k)$ is a closed set and
$\stsr\subset$ LPR$(\sum x_k)$. H. Hadwiger \cite{hadw} proved
that LPR$(\sum x_k)$ is a shifted closed additive subgroup of the
space in which the series lives. In particular this is true for
numerical series (see also \cite{sib}, exercises 3.2.2, 2.1.2 and
comments to these exercises).

By $\overline{\R}$ we denote the two point compactification of the
real line:
$$\overline{\R} = \R\cup\{-\infty,\infty\}.$$

\subsection{Main theorem for $\stsrs$}
The aim of this chapter is to prove the following result:

\begin{thm} \label{thmMain}
Let $\sum x_k\eqst a$ for the original permutation. Then $\stsr =$
LPR$(\sum x_k)$. So $\stsr$ is one of the following:
\begin{enumerate}
\item The only number $a$;

\item $\{a+\lambda\Z\}$ for some $\lambda\in\R$;

\item The whole $\R$.

\end{enumerate}

\end{thm}
\begin{proof}

Since series $\sum x_k$ converges statistically there exists a
subsequence $x_{n_k}$ such that $x_{n_k}\to 0$. From the elements of
$x_{n_k}$ we can select a subsequence $x_{n_{k_i}}$ such that
$\sum_{i=1}^\infty x_{n_{k_i}} < \infty$.

Now we can substitute all elements $x_{n_{k_i}}$ in the original
series for 0 and this will not affect the convergence since we are
subtracting an absolutely convergent series. So without loss of
generality we may assume that there are infinitely many zeros
among the original series terms.

Let us write the definition of LPR in detail:
$$
LPR(\sum x_k) = \{ x \mid \exists\,\pi\ \exists\,\{m_k\} : x =
\lim_{k\to\infty} \sum_{j=1}^{m_k} x_{\pi(j)}\}
$$
where $\pi$ is a permutation of $\N$ and $\{m_k\}$ is an increasing
sequence of indices. Let $b$ be an arbitrary element of LPR. Let
$\{m_k\}$ be a sequence from the definition corresponding to element
$b$, and such that $m_{k+1}/m_k \to\infty$. We will arrange elements
of our series in the following way:
$$
x_{\pi(1)}+\cdots+x_{\pi(m_1)} + \underbrace{0 + \cdots +
0}_{\mbox{$(m_2)^2$ times}} $$
$$+x_{\pi(m_1+1)}+\cdots+x_{\pi(m_2)} + \underbrace{0 + \cdots +
0}_{\mbox{$(m_3)^2$ times}} + \cdots
$$
We get the permutation of the series that obviously statistically
converges to $b$.

\end{proof}

\subsection{Examples}
We finish the proof by giving the examples which satisfy each case
of the theorem \ref{thmMain}.

\begin{example}
Any unconditionally convergent series in usual meaning gives us a
series with $\stsrs = \{a\}$, which corresponds case (1).
\end{example}

\begin{example}
Let the elements of series be the following:

$x_n= \left\{
  \begin{array}{ll}
    0, & \hbox{$n\neq 10^k \wedge n\neq 10^k+1$} \\
    \lambda, & \hbox{$n = 10^k$}\\
    -\lambda, & \hbox{$n = 10^k+1$}
  \end{array}
\right., k\in\N$, $n\in\N$.

Then $\stsrs = \lambda\Z$ for some $\lambda\in\R$, which corresponds
case (2).
\end{example}

\begin{example}
Any conditionally convergent series in usual meaning gives us a
series with $\stsrs = \R$, which corresponds case (4).
\end{example}

\begin{rem}
In fact the statement $\stsr =$ LPR$(\sum x_k)$ holds true for
series in any Banach space. Thus one can prove that in any
separable Banach space $\stsrs$ can be any shifted closed
subgroup.
\end{rem}

\begin{rem}
If one wants to consider $\stsrs \subset \overline{\R}$ then
modifying above argument one can prove
\begin{thm}
Let $\sum x_k\eqst a$ for the original permutation.
Then $\stsr$ is one of the following:
\begin{enumerate}
\item The only number $a$;

\item $\{a+\lambda\Z\}\cup\{-\infty, \infty\}$ for some $\lambda\in\R$;

\item The whole $\overline{\R}$;

\item The set $\{-\infty,a,\infty\}$.

\end{enumerate}

\end{thm}
\end{rem}

\section{2n-convergence}
\subsection{Introduction}
Let us say that a series $\sum_{k=1}^{\infty} x_k$ 2n-converges to
$c$ if $\lim_{n\rightarrow\infty} \sum_{k=1}^{2n} x_{k} = c$

\begin{definition}
Point $s\in\R$ belongs to the {\it 2n sum range} of the series
$\sum x_k$ if there exists a permutation $\pi: \N \to \N$ such
that $\lim_{n\rightarrow\infty} \sum_{k=1}^{2n} x_{\pi(k)} = s$.
The set of all such points is called the {\it 2n sum range} of
series $\sum x_k$ and is denoted by $\fsr$. When it is clear what
series is considered we will denote this set by $\fsrs$
\end{definition}

Consider first the following example:

\begin{example}
Series $1 + (-1) + 1 + (-1) + \cdots$

It's easy to see that this series diverges (reminder doesn't tend to
$0$). But if we consider the subsequence $S_n = \sum_{k=1}^{2n} x_k$
of it's partial sums we see that $\forall~n\in\N : S_n = 0$ and so
this subsequence converges. Notice that in order to converge
elements must go in strict pairs $1 + (-1)$ after some number of
elements. Limit of $S_n$ can be
\begin{enumerate}
  \item {$1 + 1 + 1 + (-1) + 1 + (-1) + \cdots = 2$}
  \item {$(-1) + (-1) + 1 + (-1) + 1 + (-1) + \cdots = -2$}
\end{enumerate}

It's easy to prove that
$$
\{\ S \in \R\ |\ \exists\ \pi-\mbox{permutation of}\ \N\ :\
\lim_{n\rightarrow\infty} \sum_{k=1}^{2n} x_{\pi(k)} = S \ \} =
2\,\Z
$$
\end{example}
So the statement of Riemann Rearrangement Theorem in this case of
convergence has to be modified. Surprisingly this modification and
its proof appear to be rather non-trivial and much more
complicated than in the case of statistical convergence.

Recall that $X$ is said to be {\it $\e$-separated} if all pairwise
distances between the elements of $X$ are greater than $\e$. $X$ is
said to be {\it separated} if it is $\e$-separated for some $\e
> 0$.

The aim of this chapter is to prove the following result:

\begin{thm}[Main Theorem]\label{thmMain1}
Let $\lim_{n\to\infty}\sum_{k=1}^{2n}x_k = a\in\R$. Then $\fsr$ is
one of the following:
\begin{enumerate}
\item Shifted additive subgroup of the form
$$a + \{c_1z_1 + \cdots + c_lz_l \mid z_k\in E,\ c_i \in \mathbb{Z},\ \sum_{k=1}^l
c_k \mbox{ is even}\},$$ where E of is an $\e$-separated set;

\item The whole $\R$;

\item The only number $a$.

\end{enumerate}
\end{thm}

\subsection{Reduction to a special form of the series}

We can represent the series in the following way :
\begin{equation}
\label{OriginalSeries}
x_1+(-x_1+\alpha_1)+x_3+(-x_3+\alpha_2)+x_5+(-x_5+\alpha_3)+\cdots
\end{equation}
This can be done by denoting
$$
\alpha_{k}\eqdef x_{2k-1}+x_{2k}\ (\forall~k\in\mathbb{N})
$$

Recall that the series $\sum x_k$ 2n-converges in original order,
i.e. $\lim_{n\rightarrow\infty} \sum_{k=1}^{2n} x_k = a$, so
$\sum_{i=1}^\infty \alpha_i = a$.

\begin{thm} \label{thm3.2.1}
If $\sum_{k=1}^\infty \alpha_k$ converges conditionally then
\begin{equation} \label{eqsr2}
SR_{\,2}(\sum x_k)=\R.
\end{equation}
\end{thm}

\begin{proof}
As Riemann Rearrangement Theorem says for a conditionally convergent
series $\sum_{k=1}^\infty \alpha_k$, for all $c\in \R$ there exist a
permutation of indices $\pi$ such that $\sum_{k=1}^\infty
\alpha_{\pi(k)}=c$. Consider the following arrangement of $\{x_k\}$:
$$
x_{2\pi(1)-1}+(-x_{2\pi(1)-1}+\alpha_{\pi(1)})+
x_{2\pi(2)-1}+(-x_{2\pi(2)-1}+\alpha_{\pi(2)})+\cdots
$$
It's clear that this series 2n-converges to $c$. As $c$ was
arbitrary we get (\ref{eqsr2}).
\end{proof}

\begin{definition} \label{defeq}
A series $\sum_k x_k$ is said to be equivalent to $\sum_k y_k$ if
$\sum_k |x_k - y_k| < \infty$.
\end{definition}
Remark, that if one of two equivalent series converges
(2n-converges) in some permutation then the same does the second
series and that $\fsrs(\sum_k x_k) = \fsrs(\sum_k y_k) + \sum_k
(x_k - y_k)$.

Theorem \ref{thm3.2.1}. corresponds to the case (2) of the main
theorem. Now consider what happens if $\sum_{k=1}^\infty\alpha_k$
converges unconditionally to $a$. In this case $\sum_k x_k$ is
equivalent to the following simplified series:
\begin{equation}\label{SimplifiedSeries}
x_1+(-x_1)+x_3+(-x_3)+x_5+(-x_5)+\cdots
\end{equation}

So we reduce the series ($\ref{OriginalSeries}$) to
($\ref{SimplifiedSeries}$). Changing notation we consider a series
of the form:
\begin{equation}\label{SimplifiedSeriesFinal}
x_1+x_{-1}+x_2+x_{-2}+x_3+x_{-3}+\cdots
\end{equation}
where $x_{-n} = -x_n$ and $x_n > 0$ for $n > 0$, $x_n < 0$ for $n <
0$. Denote by $X$ the set of all elements of the series (without
repetitions) and enumerate the elements of $X$ as
$$
X = \{e_i \mid i \in \mathbb{Z}\backslash\{0\}\},
$$
$e_{-n}:=-e_n$ and $e_i > 0, i\in \mathbb{N}$.  By the order of an
element $e \in X$ we mean
$$\chi(e) = \mid \{i \in \Z\backslash\{0\} \mid x_i = x \} \mid.$$

\subsection{The (basic) case of separated $X$}

\begin{lemma}\label{lemDiscrete}
Let $X$ be $\e$-separated and there are nonzero elements of infinite
order. Then
\begin{equation}\label{lemDiscreteStatement}
\begin{split}
\fsrs = \{c_1e_{j_1} + \cdots +
c_re_{j_r}~&|~\forall~k~:~\chi(e_{j_k})~=~\infty,\\&\ c_k \in\Z ,
\sum_{j=1}^r c_j~is~even\}.
\end{split}
\end{equation}
If there are no nonzero elements of infinite order then $$\fsrs =
\{0\}.$$
\end{lemma}

\begin{proof}
Denote right-hand side of (\ref{lemDiscreteStatement}) by
$\mathfrak{L}$. Let us prove that $\fsrs \subset \mathfrak{L}$.

Let $\pi\colon\N\to\Z\backslash\{0\}$ be an arbitrary bijection such
that
$$
A = \lim_{n\to\infty}\sum_{j=1}^{2n} x_{\pi(j)} \in\R
$$
Applying Cauchy Convergence Criterion to this series we get that
there exist an even number $n_0$ such that for every even $n$ and
$m$ greater or equal then $n_0$ the following inequality stands:
$|S_n -S_m| < \e$.

Consider elements of 2n partial sum sequence. Presume $n > n_0$ and
$n$ is even. We have $ |S_{n+2} - S_{n}| = |x_{\pi(n+1)} +
x_{\pi(n+2)}| < \e $ But $X$ is $\e$-separated. Since $x_j \in X$ we
get that this inequality is true if and only if $|x_{\pi(n+1)} +
x_{\pi(n+2)}| = 0$. In any other case this modulus is greater than
$\e$.
$$
(|x_{\pi(n+1)} + x_{\pi(n+2)}| = 0) \Leftrightarrow (x_{\pi(n+1)} =
- x_{\pi(n+2)})
$$
So the series has the following structure
$$
A = x_{\pi(1)} + \cdots + x_{\pi(n_0)} + x_{k_{1}} + (- x_{k_{1}}) +
x_{k_{2}} + (- x_{k_{2}}) + \cdots
$$
The first $n_0$ elements will be considered later. The other ones
come in strict pairs and if $2n$ partial sums are considered it can
be remarked that after element $S_{n_{0}}$ all of them are equal to
$S_{n_{0}}$

Now in the sum $x_{\pi(1)} + \cdots + x_{\pi(n_0)}$ consider the
elements $y_1, y_2, \dots, y_j$ of {\it finite} order. Then for
all $i\in\{1,\dots, j\}$:
$$
|\{k\in{1,\dots, n_0}\} : x_{\pi(k)}=y_i| = |\{k\in{1,\dots, n_0}\}
: x_{\pi(k)}= -y_i|
$$

That is because in the second part of the series all elements come
in strict pairs. Same number of opposite elements in finite sum
gives us zero. The only elements left are the elements with {\it
infinite} order, so $A$ has requested form
$$
A = c_1e_{j_1} + \cdots + c_re_{j_r},
$$
where $c_i \in \Z$. Moreover considering that $n_0$ was even and
although some number of pairs of elements was taken from it there
still remains an even number, we deduce that $\sum_{k=1}^l c_k$ is
even.

Why $\mathfrak{L} \subset\fsrs$?

Select a finite element $z \in \mathfrak{L}$ : $z = c_1e_{j_1} +
\cdots + c_re_{j_r}$ and create a series starting with
$$
sign(c_1)(\underbrace{e_{j_1} + \cdots + e_{j_1}}_{\mbox{$|c_1|$
times}}) + \cdots + sign(c_r)(\underbrace{e_{j_r} + \cdots +
e_{j_r}}_{\mbox{$|c_r|$ times}})
$$
and after these elements the rest of $x_k$ is settled in pairs
$$
x_{k_{1}} + (-x_{k_{1}}) + x_{k_{2}} + (-x_{k_{2}}) + \cdots
$$
It is obvious that 2n-sum of this series is $z$, and this series is
a rearrangement of ($\ref{SimplifiedSeries}$). With infinite
elements the situation is obvious. So we have managed to prove that
$\mathfrak{L} \subset\fsrs$. The second part was proved also because
the only possible real value in $\fsrs$ is $0$.

End of Lemma's proof.
\end{proof}
\subsection{Some combinatorial lemmas}
Let $M$ be a set of indices. For a bounded sequence $(x_n)_{n \in
M}$ introduce the following quantity:

\begin{equation}\label{defD(M)}
\D(M)=\D\bigl(M,(x_n)_{n \in M}\bigr)=\inf_{a \in \R}\sum_{n \in
M}|x_n - a|.
\end{equation}

If occasionally a sum consists of empty set of summands, we mean
here and below that the sum equals $0$.

\begin{lemma} \label{lemD(M)=inf}
If $\D(M) = \infty$ then one can select a finite collection of
disjoint pairs $n_k, m_k \in M$, $k = 1, 2, \ldots, s$ with
$\sum_{k=1}^s|x_{n_k} - x_{m_k}|$ being arbitrarily large.
\end{lemma}
\begin{proof}
We consider that $\{x_{n}\}_{n\in M}$ has only one limiting point
(otherwise the statement is obvious). Denote it by $a$. The fact
that $\D(M) = \infty$ implies that
$$\sum\limits_{n\in M}|x_n - a| = \infty.$$
For every $K>0$ and every $\delta>0$ there exist such $s$ that
$$\sum\limits_{k=1}^s|x_{n_k} - a| > K + \delta,$$
we just need to take a big finite part of the initial sum. Then we
can select a subsequence $\{x_{m_k}\}$  disjoint with
$\{x_{n_k}\}$ such that
$$\sum_{k=1}^s|x_{m_k} - a| < \delta,$$ this can be done since $a$
is a limit point of the sequence. This subsequences satisfy the
following inequality: $\sum_{k=1}^s|x_{n_k} - x_{m_k}|
> K$.

\end{proof}
\begin{lemma} \label{lemD(M)}
If $\D(M) < \infty$ then
\begin{enumerate}
\item\label{lemD(M)1} Either $M$ is finite, or $(x_n)_{n \in M}$
has only one limiting point.
\item\label{lemD(M)2}
The quantity $\sum_{n \in M}|x_n - a|$ in (\ref{defD(M)}) attains
its minimum in a point $a(M)$. If $M$ is infinite, then there is the
only possibility for the $a(M)$ selection: $a(M)$ must be the only
limiting point of $(x_n)_{n \in M}$. If $M$ is finite then the role
of $a(M)$ can be played by any median of $(x_n)_{n \in M}$, i.e. by
arbitrary point $a$ with the following property:
$$
\left|\{n \in M: x_n < a\}\right|=\left|\{n \in M: x_n > a\}\right|.
$$
\item \label{lemD(M)3} For every $\e > 0$ one can select a
finite collection of disjoint pairs $n_k, m_k \in M$, $k = 1, 2,
\ldots, s$ for which
\begin{equation}\label{select2)}
\sum_{k=1}^s|x_{n_k} - x_{m_k}| > \D(M)- \e.
\end{equation}
\end{enumerate}
\end{lemma}
\begin{proof}
(\ref{lemD(M)1}): Statement is obvious.

(\ref{lemD(M)2}): Assume first that $M$ is infinite. Since $\D(M)<
\infty$ there is at least one point $a$ such that $\sum_{n \in
M}|x_n - a| < \infty$. Then $x_n - a$ tends to 0 along $M$, so $a$
is the only limiting point of $(x_n)_{n \in M}$. The case of
finite $M$ is obvious.

(\ref{lemD(M)3}): First we deal with a finite $M$. In this case we
can chose $x_{n_1}$ to be the leftmost element with respect to
$a(M)$ and $x_{m_1}$ to be the rightmost, then we define $x_{n_2}$
as the leftmost element from the remaining elements $x_{m_2}$ to be
the rightmost and so on. We obtain that
$$\sum_{k=1}^s|x_{n_k} - x_{m_k}| = \D(M).$$
Suppose now $M$ to be infinite. In this case is we deal just like in
Lemma($\ref{lemD(M)=inf}$).
\end{proof}

Let $G=G_k, k \in \N$ be a disjoint collection of subsets in $\R$.
We say that $G$ is an $\e$-collection, if diameters of all the $G_k$
do not exceed $\e$. Denote $M_k = \{n \in \N: x_n \in G_k\}$;
$\D_G=\sum_{k \in \N}\D(M_k)$. Now the proof of the theorem splits
into two cases.

\subsection{Case 1 (reduction to the case of separated $X$)}
Below by distance between two sets $A,B$ of real numbers we call
$d(A,B)=\inf\{|a-b|: a \in A, b \in B\}$. Infimum of the empty set
we define to be $+ \infty$, so if at least one of $A,B$ is empty,
then $d(A,B)=+ \infty$.

\begin{lemma} \label{case1)}
Let $\{x_n\}$ have the following property: there is an $\e > 0$ such
that $\D_G < \infty$ for every $\e$-collection $G$. Then the series
$\sum x_n$ is equivalent to a series $\sum y_n$ with a separated set
of elements (like in lemma \ref{lemDiscrete}), so it satisfies the
statement of the Main theorem.
\end{lemma}
\begin{proof}
Let $\e$ satisfy the condition of the Lemma. We are going to cover
the set of values $X^+ = X \cap \R^+$ by an $\e$-collection $G$ of
intervals in such a way, that there is an $n_0$ such that for all
$n, m > n_0$ all the distances between $G_n$ and $G_m$ are bigger
than $\frac{\e}{4}$. If such $G$ is selected, put $M_k = \{n \in \N:
x_n \in G_k\}$ and denote $a_k=a(M_k) \in G_k$, $k \in \N$ the
number from Lemma \ref{lemD(M)}. In such a case the sequence $a_k$
is separated, and we can define the required symmetric sequence
$y_n, n \in \Z \setminus \{0\}$ as follows: $y_n = a_k$ for $n \in
M_k$. The set of elements of $\sum y_n$ equals $\{a_1, -a_1, a_2,
-a_2, \ldots \}$, so it is separated, and the mutual equivalence of
$\sum x_n$ and $\sum y_n$ follows from the inequality $\sum_k |x_k -
y_k| \le 2 \D_G < \infty$. So all what we need is to construct a $G$
with the property described above.

Consider covering of $X^+$ by $T_k=X^+ \cap [(k-1)\e, k\e)$ and
denote $t_k = \inf T_k$, $t^k = \sup T_k$ if $T_k \neq \emptyset$
and $t_k =  t^k = (k-1/2)\e$ if $T_k = \emptyset$. Since
$T=(T_k)_{k \in \N }$ forms an $\e$-collection, we have
$$
\sum_k |t^k - t_k | \le 2 \D_T < \infty,
$$
so in particular $|t^k - t_k | \to 0$. Select the required $n_0$ in
such a way that $|t^k - t_k | < \frac{\e}{4}$ for all $k > n_0$. For
$k \le n_0$ put $G_k = [(k-1)\e, k\e)$. Before defining $G_k$ for $k
> n_0$ let us explain the picture. We would like to take $G_k =
[t_k, t^k]$, but this can be a wrong selection, because for some $k$
both $t^k$ and $t_{k+1}$ can be very close to $k \e$ and $t_{k+1} -
t^k$ can be smaller than $\frac{\e}{4}$. But for such ``bad'' values
of $k$ the segment $[t_k, t^{k+1}]$ is of the length at most
$\frac{3 \e}{4}$, covers both the segments $[t_k, t^k]$ and
$[t_{k+1}, t^{k+1}]$, and has at least distance $\frac{\e}{4}$ from
the rest of $[t_j, t^j]$. So the required selection of $G_k$ for $k
> n_0$ can be done as follows: take all those segments $[t_j, t^j]$,
$j > n_0$,  which are far from the others (i.e. the distances to the
others are bigger than $\frac{\e}{4}$), and add all those segments
$[t_{j}, t^{j+1}]$, $j
> n_0$, where $t_{j+1} - t^j < \frac{\e}{4}$.
\end{proof}

\subsection{The remaining case}
\begin{lemma} \label{case2)}
Let $\{x_n\}$ have the the opposite to the case 1 property: for
every $\e > 0$ there is an $\e$-collection $G$ such that $\D_G =
\infty$. Then one can select a collection of disjoint pairs $n_k,
m_k \in \N$, $k = 1, 2, \ldots$ such that $|x_{n_k} - x_{m_k}| \to
0$ as $k \to \infty$ and
\begin{equation}\label{select3)}
\sum_{k=1}^\infty |x_{n_k} - x_{m_k}| = \infty.
\end{equation}
In this case $\fsr = \R$, which satisfies the statement of the Main
theorem.
\end{lemma}
\begin{proof}
For $\e = 1$ we can find an $\e$-collection $G$ such that $\D_G =
\infty$. Then applying (\ref{lemD(M)3}) of Lemma \ref{lemD(M)} one
can select a collection of disjoint pairs $n_k, m_k \in \N$, $k =
1, 2, \ldots, n_1$ such that $|x_{n_k} - x_{m_k}| < 1$ and
$\sum_{k=1}^{n_1} |x_{n_k} - x_{m_k}| > 1$. Then for $\e = 1/2$ we
select disjoint pairs $n_k, m_k \in \N$, $k = n_1 +1, \ldots, n_2$
such that $|x_{n_k} - x_{m_k}| < 1/2$ and $\sum_{k=1}^{n_2}
|x_{n_k} - x_{m_k}| > 2$. We can proceed to select the desirable
sequence of pairs. To see that $\fsr = \R$ we consider the pairs
$$(x_{n_1} - x_{m_1}), (x_{m_1} - x_{n_1}), (x_{n_2} - x_{m_2}),
(x_{m_2} - x_{n_2}), \ldots.$$ We add missing pairs of the form
$x_i - x_i$ to include all the elements into the series. Permuting
pairs (like in the Riemann rearrangement theorem) we obtain $\fsr
= \R$.
\end{proof}

\subsection{Examples}
To complete the paper we are going to demonstrate that for each of
the cases (1)\,--\,(3) of the main theorem $\ref{thmMain1}$ there
exists a series satisfying it.  To write down such examples let us
introduce more compact way of writing the series
$$
(Y,N)\eqdef \{(y_i,n_i)\mid y_i\in \mathbb{R}, n_i \in \mathbb{N}
\cup \{\infty\}\}
$$
where $n_i$ corresponds to the number of copies of $y_i$ we have in
the simplified series ($n_i$ is the order of element $y_i$) and
the following condition is satisfied
$$
y_i\not= y_j\ (\forall i \not= j).
$$
Say,
$$
\fsrs\,((1,\infty)\,,\,(-1,\infty))= \fsrs\
(1+(-1)+1+(-1)+1+(-1)+\cdots)
$$

General example : for any $\e$-separated set $E =
\{e_i\}_{i=1}^\infty$ and

$\Sg = \{(e_i, \infty)\}$ we have
$$\fsrs(\Sg) = \{c_1e_1 + \cdots + c_le_l \mid e_k\in E,
\ c_i \in \mathbb{Z},\ \sum_{k=1}^l
c_k \mbox{ is even}\}.$$ In particular consider the following
two examples.

\begin{example}
$\Sg_1 = \{(1,\infty)\,,\,(-1,\infty)\}$

Here we get that $\fsrs(\Sg_1) = \{2\Z\}$. Notice that $2\Z$ is {\it
$2$-separated}.
\end{example}

\begin{example}
$\Sg_2 =
\{(1,\infty)\,,\,(-1,\infty)\,,\,(\sqrt2,\infty)\,,\,(-\sqrt2,\infty)\}$

Applying Lemma~(\ref{lemDiscrete}) here we get that $$\fsrs(\Sg_2) =
\{a\cdot1 + b\cdot\sqrt2\},$$  where $(a + b)$ is even. It's obvious
that $\fsrs(\Sg_2)$ is dense in $\R$.
\end{example}

\begin{example}
$\Sg_3$ -- any conditionally convergent series (in the usual sense).
Obvious that it gives us $\fsrs(\Sg_3) = \R$.
\end{example}

\begin{example}
$\Sg_4 = \{\,(S_n,1)\,,\,(-S_n,1)\ |\ n\in\N\}$ , where
$S_n=\sum_{i=1}^n\frac1i$.

This series too gives us $\fsrs(\Sg_4) = \R$.
\end{example}

\begin{example} Let $\Sg_5$ be an arbitrary unconditionally
convergent series to $a\in\R$ (in usual sense). This series gives us
$\fsrs(\Sg_5) = \{a\}$.
\end{example}

\begin{rem}
If one considers convergence in $\overline{\R}$ the main theorem
must be modified as follows:
\begin{thm}
Let $\lim_{n\to\infty}\sum_{k=1}^{2n}x_k = a\in\R$. Then $\fsr$ is
one of the following:
\begin{enumerate}
\item Shifted additive subgroup of the form
$$a + \{c_1z_1 + \cdots + c_lz_l \mid z_k\in E,
\ c_i \in \mathbb{Z},\ \sum_{k=1}^l
c_k \mbox{ is even}\}\cup\{-\infty, \infty\},$$
where E of is an
$\e$-separated set;

\item The whole $\overline{\R}$;

\item The only number $a$;

\item The set $\{-\infty,a, \infty\}$.

\end{enumerate}
\end{thm}
\end{rem}

{\it Acknowledgment.} The authors are grateful to V.~Kadets for his
support, fruitful communications, and help in preparing this text. We would
also like to thank Professor Eve Oja and Professor Toivo Leiger from Tartu
University for providing us references \cite{ces} and \cite{mat}.

\end{document}